\documentclass[11pt]{article}
\usepackage{amssymb,amsmath}
\usepackage{graphicx}
\usepackage{multirow}
\newtheorem{definition}{Definition}

\newtheorem{theorem}{Theorem}
\newtheorem{lemma}{Lemma}
\newtheorem{proposition}{Proposition}

\begin{document}
\begin{center}
\large \bf A geometrically motivated parametric\\ model in manifold estimation\normalsize
\end{center}
\normalsize

\

\begin{center}
  Jos\'e R. Berrendero$^{a}$, Alejandro Cholaquidis$^{b}$,\\ Antonio Cuevas$^{a}$, Ricardo Fraiman$^{b,c}$ \\
  $^{a}$ Departamento de Matem\'aticas, Universidad Aut\'onoma de Madrid, Spain\\
  $^{b}$ Centro de Matem\'atica, Universidad de la Rep\'ublica, Uruguay\\
  $^{c}$ Departamento de Matem\'atica, Universidad de San Andr\'es, Argentina\\
\end{center}

\section*{Abstract}  {\small The general aim of manifold estimation is reconstructing, by statistical methods, an $m$-dimensional
compact manifold $S$ on ${\mathbb R}^d$ (with $m\leq d$) or estimating some relevant quantities related to the geometric properties of $S$.
We will assume that
the sample data are given by the distances to the $(d-1)$-dimensional manifold $S$ from points randomly chosen
on a band surrounding $S$, with $d=2$ and $d=3$.
The point in this paper is to show that, if $S$ belongs to a wide class of compact
 sets (which we call \it sets with polynomial volume\rm), the proposed
 statistical model leads to a relatively simple parametric formulation. In this setup, standard
 methodologies (method of moments, maximum likelihood) can be used to estimate some interesting
 geometric parameters, including curvatures and Euler characteristic. We will particularly
 focus on the estimation of the $(d-1)$-dimensional boundary measure (in Minkowski's sense) of $S$.

It turns out, however, that the estimation problem is not straightforward since the standard estimators show a remarkably pathological behavior: while they are consistent and asymptotically normal, their expectations are infinite. The theoretical and practical consequences of this fact are discussed in some detail.}

\

\noindent{\footnotesize \bf AMS 2010 Subject Classification\rm: 62F10, 62H35.}

\noindent{\footnotesize \bf Key words and phrases\rm: Estimation of boundary length; estimation of curvature; distance to boundary; volume function; remote sensing.}

\noindent{\footnotesize This work has been partially supported by Spanish Grants MTM2010-17366 (Authors 1, 3 and 4) and CCG10-UAM/ESP-5494 (Authors 1 and 3).}

\newpage

\section{Introduction}

\

\noindent \it Some background: manifold estimation\rm

Let $S$ be an $m$-dimensional
compact manifold in ${\mathbb R}^d$, that is, a compact subset of ${\mathbb R}^d$ in
which every point has a neighborhood which is homeomorphic to an open $m$-ball in ${\mathbb R}^m$, where $1\leq m\leq d$.

For $r>0$ define the \it $r$-parallel\/ \rm set (or \it $r$-offset\rm),
$$
B(S,r)=\{x: {\mathcal D}(x,S)\leq r\}=\bigcup_{x\in S}B(x,r),
$$
where ${\mathcal D}(x,S)=\mbox{inf}_{y\in S}\Vert x-y\Vert$ and
$B(x,r)$ denotes the Euclidean closed ball of center $x$ and radius $r$.

A natural approach to tackle the estimation of $S$ is to assume that the sample points
$X_1,\ldots,X_n$ are randomly drawn ``around'' $S$. More formally, these points
could arise as noisy versions of observations randomly chosen on $S$, that is $X_i=\xi_i+Z_i$, where the $\xi_i$ are random points chosen on the boundary $\partial S$ and $Z_i$ are iid random observations from a noise random variable $Z$. This is the \it additive noise model\rm, used (under different assumptions for the noise variables $Z_i$)  by Genovese \it et al\rm. (2012a, 2012b) and Niyogi \it et al\rm.
(2008).

We will consider here a slightly different \it clutter noise model\/ \rm which, in his more general formulation [Genovese \it et al\rm. (2012b)], assumes that the sample observations $X_1,\ldots,X_n$ come from a distribution $(1-\alpha) F+\alpha G$, where $F$ is supported on $S$, $G$ is uniformly distributed on a compact set ${\mathcal K}\subset {\mathbb R}^d$ and $\alpha\in [0,1]$. To be more specific, be will consider the case of ``extreme'' noise contamination where $\alpha=1$ and the noise support ${\mathcal K}$ is the topological closure
of an outside band $B(S,R)\setminus S$ surrounding $S$.

The problem of estimating $S$ under such sample models is a relatively new subject, of increasing interest, usually called \it manifold estimation\rm, which can be included in the broader field of \it manifold learning\rm; see Dey (2007) for a recent general reference.

Manifold estimation is closely related, in the statement and methodology of the problem,
to the theory of \it set estimation\/ \rm (see Cuevas and Fraiman (2009) for a recent survey)
and, more specifically, to \it boundary estimation\rm: see Cuevas and Rodr\'{\i}guez-Casal (2004). There,
the problem is essentially to estimate the boundary $\partial S$ of a given set $S\subset{\mathbb R}^d$ from an
iid sample of a probability distribution with support $S$. As it sometimes happens in the statistical research,
the focus of boundary (or set) estimation soon moved from the primary target of estimating the boundary (or set)
 itself, to other related goals that can be formalized in terms of estimation of appropriate functionals.
 A relevant example is the estimation of the $(d-1)$-dimensional measure of $\partial S$. See Cuevas et al. (2007), Pateiro-L\'opez
 and Rodr\'{\i}guez-Casal (2008), Armend\'ariz et al. (2009) and Jim\'enez and Yukich (2011). In all these references,
 the sample model is somewhat different from the original simple iid situation mentioned above since  the available sample information consists of random points drawn inside and outside $S$. In some sense, the present paper goes along similar lines in the problem of manifold estimation as our main concern here will be the estimation of the $(d-1)$-dimensional measure of a manifold $S$ with dimension $(d-1)$.

\

\noindent \it The manifold problem and the solid problem. Our sampling model(s)\rm

Given a set $A\subset {\mathbb R}^d$, we define the $(d-1)$-dimensional \it Minkowski content\/ \rm of its topological boundary $\partial A$ by
\begin{equation}
L(\partial A)=\lim_{\epsilon\to 0}\frac{\mu(B(\partial A,\epsilon))}{2\epsilon},\label{Minkowski}
\end{equation}
provided that this limit is finite; here $\mu$ denotes the $d$-dimensional Lebesgue measure.

Likewise, the one-sided (outer) Minkowski content of $A$ is defined by
\begin{equation}
L^+(A)=\lim_{\epsilon\to 0}\frac{\mu(B(A,\epsilon)\setminus A)}{\epsilon}.\label{outer-Minkowski}
\end{equation}
Assuming that $A$ has a Lipschitz boundary, it can be proved (see Ambrosio et al., 2008, Theorem 5, for a precise statement) that $L^+(A)=L(\partial A)$.

In this paper, we will consider two slightly different problems
whose statistical treatment turns out to be essentially identical.
First, the estimation of $L(\partial S)$ when $S$ is a $(d-1)$-dimensional,
smooth enough, compact manifold (so $\mu(S)=0$ and $\partial S=S$), will be called
the \it manifold problem\/ \rm or the \it manifold model\rm. Second, the estimation of $L^+(\partial S)$, when $S$ is
a $d$-dimensional set (with non-empty interior), will be called the \it solid
problem\/ \rm or the \it solid model\rm.

As an important difference with respect to the general manifold estimation problem mentioned above, we will assume (in both models) that the
sample data consist of the
distances $D_i={\mathcal D}(X_i,S)$ $i=1,\ldots,n$ to $S$ from points $X_1,\ldots,X_n$ uniformly drawn on the parallel set $B(S,R)$ but \it outside\/ \rm $S$. Therefore, whereas in the manifold problem (where typically $\mu(S)=0$) this amounts to draw the random points $X_i$
on the whole parallel set $B(S,R)$, in the solid problem (where $\mu(S)>0$) we will assume that the $X_i$ are drawn on $B(S,R)\setminus S$.
This distinction makes sense as, in practice, it is reasonable to assume that we are just allowed to observe $S$ ``from the outside''. As we will see, the mathematical treatment is essentially identical in both cases, with just a few minor differences. For this reason we will denote, with some abuse of notation,  $L_0$ in both cases the value of the target parameter. See expressions (\ref{pv-2}) and (\ref{pv-3}) below for details.

These sampling models can be motivated in terms of \it remote sensing\rm: we could think that we are able to measure (with a sonar device,
for example) the distance from the outside points $X_i$ to the surface $S$ or to the solid $S$.

\

\noindent \it The structure of this paper\rm.

Some necessary concepts, related to the structure of the volume function and its geometric and analytic interpretations are reviewed in Section 2.

The basic geometric assumption, as well as the main theoretical results concerning the estimation of the boundary measure $L_0$, are established in Section 3. To be more specific, we show that, according to the proposed model, the distribution of the random variable $D={\mathcal D}(X,S)$ belongs to a parametric family indexed by $L_0$ in such a way that the theoretical expressions  for the asymptotic distributions of both the maximum likelihood and the moment estimator of $L_0$ can be explicitly obtained. We consider the two-dimensional case $d=2$ (where $S$ is a curve in the manifold model and a planar domain in the solid model) and the three-dimensional one $d=3$ (where $S$ is a surface or solid body, respectively). In the case $d=3$, besides the estimation of $L_0$, we can  also tackle the estimation of a parameter, denoted by $M$, which can be interpreted as the \it integrated mean curvature\/ \rm of $S$.

In Section 4 we show that the standard estimators are in some sense, pathological. In particular, the moments estimator (in spite of being consistent and asymptotically normal) has an infinite expectation. This entails
that the usual mean square error is no longer a suitable criterion to measure the performance of these estimators. Hence an alternative error criterion is proposed. Also other
estimation methods, aimed to overcome the infinite expectation pathology are considered.

Section 5 is devoted to a small simulation study.

Section 6 includes some discussion and a few final remarks.


\section{Some geometric preliminaries. The volume function}

In what follows the \it volume function\/ \rm
$V(r)=V(r;S)=\mu(B(S,r))$ plays an outstanding role. It appears in
a natural way in different topics related to stochastic geometry
and geometric measure theory; see, e.g., Hug \it et al\rm. (2004),
Ambrosio \it et al\rm. (2008) and Villa (2009) for recent references. In
set estimation $V(r)$ arises also as an auxiliary tool to obtain
convergence rates with respect to the Hausdorff metric; see, e.g.,
Walther (1997). An additional statistical application of $V(r)$
will be presented in this paper.

The discussion below involves the use of some classical, though non-trivial, concepts from differential geometry
and geometric measure theory. This section is devoted to briefly outline them. We just introduce the main results and concepts,
pointing out their intuitive meanings, and refer to some standard references for additional details.

\

\noindent \it The Steiner formula\rm

The systematic study of the volume function $V(r;S)$ goes back to the nineteenth century. The best known result about this function is maybe the classical Steiner's  (1840) formula  whose $d$-dimensional version is as follows:  If $S\subset {\mathbb R}^d$ is a compact convex set,
 then the corresponding volume function is a polynomial in $r$ of degree $d$,
\begin{equation}
V(r;S)=\mu(B(S,r))=\sum_{j=0}^d r^{d-j}\omega_{d-j}V_j(S),\label{Steiner}
\end{equation}
where $\omega_k$ denotes the ($k$-dimensional) volume of the Euclidean unit ball in ${\mathbb R}^k$ (with $\omega_0=1$) and the
coefficients $V_0(S)\ldots,V_d(S)$ are the so-called ``intrinsic volumes'' of $S$. In particular, $V_d(S)=\mu(S)$
is the volume of $S$, $V_0(S)=1$ and $V_{d-1}(S)=L^+(S)/2$.

In the cases $d=2$ and $d=3$, we will express the Steiner formula with the notations,
\begin{equation}
V(r;S)=\mu(S)+L_0r+\pi r^2,\label{Steiner2}
\end{equation}
and
\begin{equation}
V(r;S)=\mu(S)+L_0r+M r^2+4\pi r^3/3,\label{Steiner3}
\end{equation}
respectively. The value $M$
in (\ref{Steiner3}) coincides with the ``integrated mean curvature'' of $\partial S$.

\

\noindent \it Sets of positive reach. Federer's volume formula and its geometrical interpretation\rm

The appearance of a measure of curvature in (\ref{Steiner3}) is not by chance. This point was clarified by Federer (1959) in a celebrated paper which, in many respects, can be considered as the pioneering reference in Geometric Measure Theory. In that paper, a generalization of the Steiner formula (\ref{Steiner3}), together with a deep interpretation of the corresponding polynomial coefficients, was established.

Federer's result is valid for a broad class of sets having a \it positive reach\/ \rm
 property. This is a fairly intuitive smoothness condition which does not involve any explicit differentiability assumption. The \it reach\/ \rm of a (closed) set $S$, $\mbox{reach}(S)$, is defined as the largest $r$ (possibly $\infty$) such that if ${\mathcal D}(x,S)<r$ then $S$ contains a unique point nearest to $x$. If $r_0=\mbox{reach}(S)>0$ then $S$ is said to have positive reach.

As a combination of Theorems 5.6 and 5.19 in Federer (1959) we have the following clean and powerful result:

\

\noindent \sc Federer's Theorem.- \it If $S\subset{\mathbb R}^d$ is a compact set with $r_0=\mbox{reach}(S)>0$, then there exist
unique $\Phi_0(S),\ldots,\Phi_d(S)$ such that
\begin{equation}
V(r;S)=\sum_{j=0}^d r^{d-j}\omega_{d-j}\Phi_j(S),\ \mbox{for } 0\leq r<r_0.\label{Federer}
\end{equation}
Moreover, $\Phi_0(S)$ coincides with the so-called Euler characteristic of $S$ (see below) which, in particular, is a topological invariant\rm.

\

It is readily seen that $\Phi_d(S)=\mu(S)$ and $\Phi_{d-1}(S)=L^+(S)/2$, where $L^+$ is defined in (\ref{outer-Minkowski}). The meaning of the remaining coefficients $\Phi_j(S)$ is also carefully addressed in Federer (1959) by showing that they can be interpreted as the \it total curvatures\/ \rm of $S$.

 The above theorem is a considerable extension of the Steiner formula. Note that it applies of course to any convex compact set since a closed set $S$ is convex if and only if $\mbox{reach}(S)=\infty$. Moreover, as Federer (1959, Section 4) points out the class of sets with positive reach \it ``contains (...) all those sets which can be defined locally by means of finitely many equations, $f(x)=0$, and inequalities, $f(x)\leq 0$, using real valued continuously differentiable functions, $f$, whose gradients are Lipschitzian and satisfy a certain independence condition''\rm. Note that if $S$ has a positive reach condition then $\partial S$ can have ``outward peaks'' but the ``inward'' (non-differentiable) peaks are ruled out.

\

\noindent \it The Euler characteristic\rm

As indicated above, the total curvature $\Phi_0(S)$ in (\ref{Federer}) equals the \it Euler characteristic of $S$\rm.
This is an important, integer-valued, quantity which provides useful information on some geometric aspects of a surface or, more in general, of a topological space.

Let us recall that a  \it Riemannian manifold\/ \rm is just a
differentiable manifold in which every tangent space is equipped
with an inner product with an associated  Riemannian metric
which varies smoothly from point to point.

The formal relation of the Euler characteristic  with the notion of curvature is given by the Gauss-Bonnet theorem.
The simplest version of this result states that the total Gaussian curvature of a  compact two-dimensional Riemannian manifold $T$  without boundary  is equal to $2\pi\chi(T)$ where $\chi(T)$ denotes the Euler characteristic of the surface.
The result can be extended to even-dimensional manifolds (the Euler characteristic of an odd dimensional
compact manifold is zero).
This is a striking fact since, in principle, the curvature is a notion that depends on local properties of the surface (relying on differentiability properties) and the Euler characteristic is a global, topological invariant, which means that it does not change by bijective bi-continuous transformations.

Let us now briefly recall some basic facts about the Euler characteristic. A more complete discussion can be found in the book by Hatcher (2002).
The simplest definition of $\chi$ can be given for polyhedral surfaces $T$ in ${\mathbb R}^3$. In this case  $\chi(T)=\mbox{(number of vertices)}-\mbox{(number of edges)}+\mbox{(number of faces)}$. It is a well-known classical result that if $T$ is the boundary of a convex polyhedron then
$\chi(T)=2$. The Euler characteristic can be defined for any subset of ${\mathbb R}^d$ in such a way that it is a topological invariant. As a consequence of this invariance we also have that, for the two-dimensional sphere $S^2$ (i.e. the boundary of the three-dimensional ball)  $\chi(S^2)=2$ and the same holds for any compact orientable surface homeomorphic to the sphere.

In the $d$-dimensional case we have that $\chi(S^d)=1+(-1)^d$ so that it is always 0 or 2.

The Euler characteristic is also related to other invariants. For example, for connected orientable compact surfaces without boundary, we have $\chi=2-2g$, where $g$ is the \it genus\/ \rm of the surface, which intuitively coincides with the number of ``handles''. Thus, $\chi=0$ for a torus, $\chi=-2$ for a double torus (with two handles) and so on.

The value of Euler's characteristic is also explicitly known for many other interesting sets in ${\mathbb R}^d$, not necessarily curves or surfaces. For example, it is known that if $S\subset{\mathbb R}^d$ is a  ``solid'' ball, then $\chi(S)=1$. In fact the same is true for any contractible set (i.e., homotopy equivalent to a point). It follows from the previous discussion that the class of compact sets $S$ in ${\mathbb R}^d$ with $\chi(S)=1$ is extremely wide.

\

\section{Statistical results: parametric estimation of some geometric quantities}\label{stat}

\

\noindent \it The statistical interpretation of the volume function\rm

According to the statistical model(s) established in the introduction, we will always assume that our sample data consist of iid observations $D_1,\ldots,D_n$ from the distance variable $D={\mathcal D}(X,S)$ so that $D_i={\mathcal D}(X_i,S)$ where $X_1,\ldots,X_n$ are iid random variables with uniform distribution on the band $B(S,R)\setminus S$. We will simultaneously consider the \it manifold model \/ \rm where $S$ will be a $(d-1)$-dimensional manifold with $\mu(S)=0$ and the \it solid model\/ \rm where $\mu(S)>0$. In both cases the main target will be to estimate the surface measure $L_0$.

 The following proposition is just a reformulation, in terms of our statistical model, of some results proved by Stach\'o (1976). It is included here for the sake of completeness.

 \begin{proposition}\label{stacho}  Let $S\subset{\mathbb R}^d$ be a compact set and $R>0$ a fixed constant. Given a random variable $X$ uniformly distributed on the band $B(S,R)\setminus S$, define the (Euclidean) distance variable $D={\mathcal D}(X,S)$. Denote by $F(r)={\mathbb P}(D\leq r)$, for $0\leq r\leq R$, the distribution function of $D$.

\begin{itemize}
\item[(a)]  The distribution function $F$ is given by
\begin{equation}
F(r)=\frac{V(r)-\mu(S)}{V(R)-\mu(S)},\ 0\leq r\leq R,\label{FV}
\end{equation}
where $V(r)=\mu(B(S,r))$ is the volume function associated with $S$. Moreover, $F$ is absolutely continuous and differentiable except for, at most, a countable set of points. In particular, it can be expressed as the integral of its derivative, $F(r)=\int_0^r F^\prime(t)dt$ where $F^\prime(r):=f(r)=\frac{V^\prime(r)}{V(R)-\mu(S)}{\mathbb I}_{[0,R]}(r)$ a.e. ($\mu$) is the density function of $D$.

\item[(b)] For every $r>0$ the left and right hand-side derivatives $F_-^\prime$ and $F_+^\prime$ do exist. Moreover they are continuous from the left and from the right, respectively and fulfill $F_-^\prime\geq F_+^\prime$.
\item[(c)] For all $r>0$ there exists the Minkowski measure $L(\partial B(S,r))$ and
\begin{equation}
L(\partial B(S,r))=\frac{V(R)-\mu(S)}{2}\left(F_-^\prime(r)+F_+^\prime(r)\right)\ \mbox{for\ } 0<r<R\label{Stacho}
\end{equation}

\end{itemize}
\end{proposition}

\rm

\noindent
\sc Proof\rm: (a) and (b) Since $X$ is uniformly distributed, we have
{\small \begin{equation*}
F(r)={\mathbb P}\{D\leq r\}={\mathbb P}\{X\in B(S,r)\setminus S\}=\frac{\mu(B(S,r)\setminus S)}{\mu(B(S,R)\setminus S)}=\frac{V(r)-\mu(S)}{V(R)-\mu(S)}.
\end{equation*}}
Now,  statements  (a) and (b) concerning the absolute continuity and differentiability properties of $V$
follow directly from Lemma 2 in Stach\'o (1976). In fact, these properties are established in general for the  so-called \it functions of Kneser type\rm, and it is shown that the volume function belongs to that class. This means that  $V(\lambda b)-V(\lambda a)\leq \lambda^d(V(b)-V(a))$, for all $0\leq a\leq b$, $\lambda\geq 1$.

Result (c) is just Theorem 2 in Stach\'o (1976) rewritten in our statistical framework.
\hfill $\square$

\

\noindent \it The basic geometric assumption: sets with polynomial volume\rm

According to Proposition \ref{stacho}, the simpler the structure of $V(r)$ the easier the statistical problem stated in the introduction.
The discussion in the previous section suggests that, concerning $V(r)$,  we cannot expect anything simpler than the polynomial structure given by Steiner's theorem. However, as we have also pointed out, there is no need to assume that $S$ is convex in order to get a polynomial volume function (at least on a given interval).

This lead us in a natural way to the following definition.

\begin{definition}\label{PV}
We will say that $S\subset{\mathbb R}^d$ is a set of polynomial volume, of type 1, on the interval $[0,R]$ if the
volume function has an expression of type
\begin{equation}
V(r;S)=\omega_d r^d+\sum_{j=1}^d r^{d-j}\omega_{d-j}\Phi_j(S),\ \mbox{for } 0\leq r<R,\label{PVR}
\end{equation}
where $\omega_k$ denotes the volume of the unit ball in ${\mathbb R}^k$ and $\Phi_j(S)$ are appropriate coefficients. The family of sets in ${\mathbb R}^d$ fulfilling this property will be denoted by ${\mathcal P}{\mathcal V}(R,d)$. More generally, we could also define the class ${\mathcal P}{\mathcal V}(R,d,\Phi_0)$ of sets $S\subset{\mathbb R}^d$ of polynomial volume, of type $\Phi_0$, by imposing that their volume functions have an expression such as (\ref{PVR}) where the term $\omega_d r^d$ is replaced with $\Phi_0\omega_d r^d$.
\end{definition}

As a consequence of Steiner's theorem, the class ${\mathcal P}{\mathcal V}(R,d)$ includes that of compact convex sets in ${\mathbb R}^d$ but, from Federer's theorem, it also includes the much broader class of compact sets with reach $R$ and Euler characteristic 1.

Since the class of sets with positive reach is by far the best known
class of sets with a polynomial volume on a interval, it is
natural to ask whether there exists a simple characterization of those sets that having a polynomial volume
but still do not fulfill the positive reach property. As far as we
know, this is still an open question (see Heveling et al., 2004
for interesting closely related issues). It is easy to construct simple examples of such sets. Thus, the polygonal $S$ joining the points $(-1,1)$, $(0,0)$ and $(1,1)$ belongs to the family ${\mathcal P}{\mathcal V}(R,2,\Phi_0))$ with $R<1$ and $\Phi_0=\frac{5}{4}-\frac{1}{\pi}$. The same holds for the non-convex pentagon $[0,1]^2\setminus T$, where $T$ denotes the (open) triangle whose vertices are $(0,1)$, $(\frac{1}{2},\frac{1}{2})$ and $(1,1)$.

The set (b) in Figure 1 is defined as the unit circle in ${\mathbb R}^2$ minus the cone with center $(0,0)$ and angle $\rho\in(0,\pi/2)$. A direct calculation shows that in this case the volume function is 
$$
V(r)=\left(\frac{3\pi-\rho}{2}-\frac{1}{\tan{(\rho/2)}}\right)r^2+(2\pi-\rho+2)r+\pi-\frac{\rho}{2},\ \mbox{for } 0\leq r\leq\tan{(\rho/2)}.
$$

Heveling et al. (2004)
present a general construction of non-convex sets in
$\mathbb{R}^d$ ($d\geq 3$) with reach equal to zero and with
polynomial volume function $V(r)$ for any $r\geq 0$. Examples of them
are those presented in Figure 1, (c) and (d). The first one is just 
the union of two touching balls, $S=B\big((0,0,1),1\big)\cup B\big((0,0,-1),1\big)$. It can be seen that
\begin{eqnarray*}
V(r) &=& \frac{8}{3}\pi(1+r)^3- \Big(\frac{4}{3}\pi(1+r)^3-2\pi(1+r)^2+\frac{2\pi}{3}\Big)\\ 
&=&\frac{4}{3}\pi(1+r)^3+2\pi(1+r)^2-\frac{2\pi}{3}
= \frac{4}{3}\pi r^3 +6\pi r^2+8\pi r+\frac{8}{3}\pi
\end{eqnarray*}
The set (d) in Figure 1 can be defined as $S=B(A,1)$, where $A$ is the union of the closed segment joining the points $(0,0,-1/2)$ and $(0,0,-1)$ with the point $(0,0,1)$. It can be proved that in this case
$$
V(r)=\frac{4}{3}\pi(1+r)^3+\frac{3}{2}\pi(1+r)^2.
$$
As a conclusion, the cases (a) and (b) in Figure 1 provide examples of sets in ${\mathbb R}^2$ with polynomial volume 
but not of type 1, that is they are not in ${\mathcal P}{\mathcal V}(R,2)$ since the value $\Phi_0$ in the highest order term 
$\Phi_0\omega_d r^d$ of the polynomial volume function is not 1. On the other hand, the cases (c)  and (d) correspond to sets with reach 0 but belonging to
${\mathcal P}{\mathcal V}(R,2)$ for all $R>0$. 

\

\begin{figure}[h]
\includegraphics
[scale=.6]{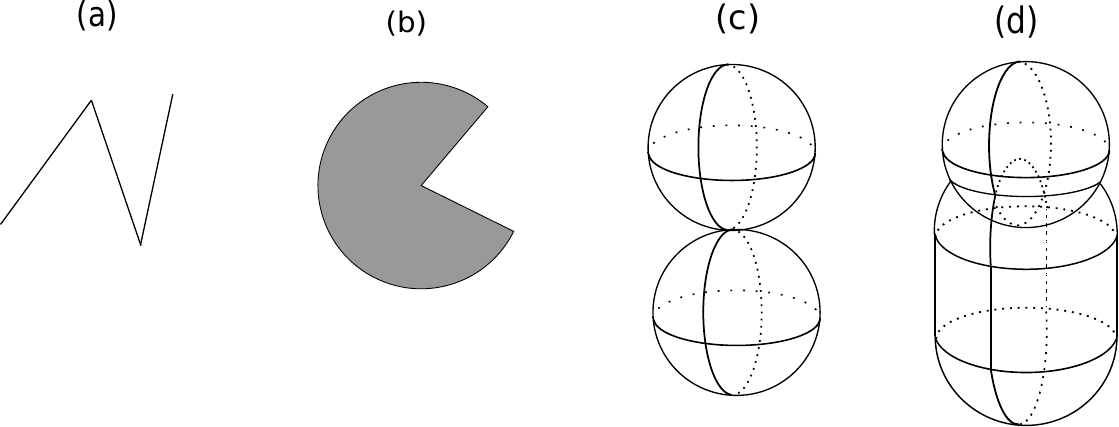}
\caption{(a) y (b) are planar sets not in $\mathcal{PV}(R,2)$  whose volume function is polynomial.
(c) and (d) are sets with reach 0,
in ${\mathcal P}{\mathcal V}(R,3)$ for all $R\geq 0$. }
\end{figure}

Throughout the rest of the paper we shall concentrate on the cases $d=2$ and $d=3$ (though the basic ideas can potentially be extended to general dimensions). So we will deal with the classes ${\mathcal P}{\mathcal V}(R,2)$ and ${\mathcal P}{\mathcal V}(R,3)$ for which the expressions of $V(r)$ on the interval $[0,R]$ are of type (\ref{Steiner2}) and (\ref{Steiner3}), respectively. Let us recall that these classes include all sets with positive reach and Euler's characteristic 1. The more general cases ${\mathcal P}{\mathcal V}(R,2,\Phi_0)$ and ${\mathcal P}{\mathcal V}(R,3,\Phi_0)$ can be handled in a  similar way, just incorporating $\Phi_0$ as an additional parameter in the estimation procedure (in case it were not known in advance).

\

\noindent\it The two-dimensional case\rm

Let us first consider the case where $S\subset {\mathbb R}^2$. In this case, our only estimation target is $L_0$. The following result provides two alternative expressions for the distribution of the
random variable ``distance to the boundary of $S$'', $D$, defined above. We assume that $S$ belongs to the class ${\mathcal P}{\mathcal V}(R,2)$ of sets with polynomial volume given by
\begin{equation}
V(r;S)=\mu(S)+L_0r+\pi r^2,\ 0\leq r<R.\label{pv-2}
\end{equation}

\begin{proposition}
\label{proposition:pv2}

The random variable $D$  is absolutely continuous with density function
\begin{equation}
f(r)=\frac{L_0 +2\pi r}{L_0R+\pi R^2},\ 0\leq r<R.\label{param1}
\end{equation}
An alternative expression for this density  is
\begin{equation}
f(r) = \lambda f_1(r) + (1-\lambda) f_2(r),\ \ \ 0\leq r<R,\label{mixture}
\end{equation}
 where $\lambda = L_0 (L_0 +\pi R)^{-1}$, $f_1$ is  the density function of a random variable $Z_1$, uniform in $(0,R)$ and $f_2$ is the density function of $RZ_2$,  where $Z_2$ follows a Beta distribution with parameters $\alpha=2$ and $\beta=1$.

\end{proposition}

\rm
\noindent
\sc Proof\rm: Expression (\ref{param1}) follows directly from Proposition \ref{stacho}. Expression (\ref{mixture}) is just a simple reformulation of (\ref{param1}). \hfill $\square$

\

In order to gain some insight on the geometric meaning of (\ref{mixture}), let us consider the simple case of a square $S=[0,1]^2$.
While the distance $D$ from  those ``regular'' points in $B(S,R)\setminus S$ not projecting on any of the vertices of $S$
follows a uniform distribution $f_1$, the density $f_2$ accounts for the remaining points whose projection is one vertex. For more complicated sets one could think that (\ref{mixture}) reflects the mixture between ``flatness'' (the $f_1$ term)  and ``curvature''(the $f_2$ term)  in the boundary of $S$.

We are now ready to consider the estimation of $L_0$. Let us first analyze the solution provided by the classical method of moments. The following theorem shows that, at first sight, this procedure works reasonably well, in the sense that the expression of the estimator is not too complicated and the asymptotic distribution is easy to identify. However, as we will see in the next section, a rather surprising property comes up.

\begin{theorem}
\label{theorem:2dim}
Under the assumption (\ref{pv-2}), we have that the estimator of $L_0$ by the method of moments from a sample $D_1,\ldots,D_n$ of $D$ is given by
\begin{equation}
\tilde{L}_0 =\frac{2\pi R}{3}\, \frac{2R-3\bar{D}}{2\bar{D}-R},\label{emm}
\end{equation}
where $\bar D$ denotes the sample mean of $D_1,\ldots,D_n$.

This estimator is asymptotically normal. More precisely, we have
\begin{equation}
\sqrt{n}(\tilde{L}_0-L_0)\stackrel{w} \longrightarrow\mbox{N} (0,\sigma^2_{\tilde{L}_0}),\label{emman}
\end{equation}
where $\stackrel{w}\longrightarrow$ stands for convergence in law and
\begin{equation}
\sigma^2_{\tilde{L}_0} = (L_0+\pi R)^2\left[3\left(1+\frac{L_0}{\pi R} \right)^2-1\right].\label{emmav}
\end{equation}
\end{theorem}

\noindent \sc Proof\rm:
(a) First, we compute the expected distance:
\[
\mu := {\mathbb E}(D) = \int_0^R r \frac{L_0 + 2\pi r}{L_0 R + \pi R^2} dr = \frac{3L_0R+4\pi R^2}{6(L_0+\pi R)}
\]
The moment estimator, $\tilde{L}_0$ is defined to be the solution in ${L}_0$ of the equation
\begin{equation}
\bar{D} = \frac{3{L}_0R+4\pi R^2}{6({L}_0+\pi R)},\label{barD}
\end{equation}
where $\bar D$ denotes the sample mean of the sample $D_1,\ldots,D_n$. Thus,  from (\ref{barD}) we easily get (\ref{emm}).

From the Central Limit Theorem applied to $\bar D$ we have,
\[
\sqrt{n}(\bar{D} - \mu) \stackrel{w}\longrightarrow\mbox{N} (0,\mbox{Var}(D))
\]
where, after some algebra, it is not difficult to show that
\[
\mbox{Var}(D) = \frac{R^2(3L_0^2 + 6\pi RL_0 + 2\pi^2 R^2)}{36(L_0 + \pi R)^2}.
\]
Now observe that $\sqrt{n}(\tilde{L_0}-L_0) = \sqrt{n}[g(\bar{D}) - g(\mu)]$, where
$g(u)= \frac{2\pi R}{3}\, \frac{2R-3u}{2u-R}$. It is easy to check that
\[
\mu = \frac{3L_0 R + 4\pi R^2}{6(L_0 + \pi R)} = \frac{R}{2} + \frac{\pi R^2}{6(L_0+\pi R)} > \frac{R}{2}.
\]
Thus, $\mu\neq R/2$ and $g'(\mu)\neq 0$. Notice also that
\begin{equation*}
\label{eq.2mu-R}
\mu - R/2 = \frac{\pi R^2}{6(L_0 + \pi R)}:=\delta_0>0.
\end{equation*}
Therefore, using the standard delta-method for $g$ restricted to the interval $[\mu-\delta_0/2,\mu+\delta_0/2]$, [e.g. Lehmann and Casella (1998), Th. 8.12, p. 58] we conclude
$$
\sqrt{n}(\tilde{L}_0-L_0)\stackrel{w} \longrightarrow\mbox{N} (0,(g^{\prime }(\mu))^2\mbox{Var}(D))
$$
which leads to (\ref{emman}). \hfill $\square$

\

The next theorem is devoted to analyze the properties of the maximum likelihood estimator $\widehat{L}_0$. Unlike the moment estimator,
$\widehat{L}_0$ has no explicit expression but, as we will see, it is slightly more efficient.

\begin{theorem}
\label{theorem:2dim-mv}
Under the assumption (\ref{pv-2}), we have
 that the maximum likelihood estimator of $L_0$, $\widehat{L}_0$,  appears as the solution of the likelihood equation
\begin{equation}
\frac{1}{n}\sum_{i=1}^n \left(\frac{1}{\widehat{L}_0 + 2\pi D_i}\right)
= \frac{1}{\widehat{L}_0 + \pi R}.\label{le}
\end{equation}
 This estimator is asymptotically normal, that is,
\begin{equation}
\sqrt{n}(\widehat{L}_0-L_0) \stackrel{w}\longrightarrow\mbox{N} (0,\sigma^2_{\widehat{L}_0}),\label{leasymp}
\end{equation}
where
\[
\sigma^2_{\widehat{L}_0} = (L_0 + \pi R)\left[\frac{1}{2\pi R}\log\left(1+\frac{2\pi R}{L_0}\right) - \frac{1}{L_0 + \pi R} \right]^{-1}
\]
coincides with the Fr\'echet-Cramer-Rao bound (given by the inverse
of Fisher's information measure).
\end{theorem}

\noindent \sc Proof\rm:

The likelihood equation (\ref{le}) follows directly by calculating the derivative with respect to $L_0$ of the log-likelihood,
\[
\ell(L_0;D_1,\ldots,D_n) = -n\log(L_0 R + \pi R^2) + \sum_{i=1}^n \log(L_0 + 2\pi D_i).
\]
As for (\ref{leasymp}), we will use the standard result on asymptotic normality of the maximum likelihood estimation which can be found in many standard textbooks. We will use the version given in Lehmann and Casella (1998), Th. 3.10, p. 449.  According to this result, a conclusion of type (\ref{leasymp}) can be obtained, for a general one-parameter family given by the (Lebesgue) densities $f(x;\theta)$, $\theta\in \Omega$, under the following regularity conditions:
\begin{itemize}
\item[(i)] The parameter space $\Omega$ is an open interval (not necessarily finite).
\item[(ii)] The support of the distributions in the parametric family does not depend on $\theta$, so that the set $A=\{x:f(x;\theta)>0\}$ is independent of $\theta$
\item[(iii)] For every $x\in A$ the density $f (x;\theta)$ is three times differentiable with respect to $\theta$,
and the third derivative is continuous in $\theta$.
\item[(iv)] The integral
$\int f (x;\theta)dx$ can be three times differentiated under the integral sign.
\item[(v)] The \it Fisher information\/ \rm $I(\theta)={\mathbb E}_\theta\left[(\frac{\partial}{\partial \theta}\log f(X;\theta))^2\right]$ fulfills $0<I(\theta)<\infty$.
\item[(vi)] For any given $\theta_0\in\Omega$, there exists a positive number $c$ and a function $M(x)$
(both of which may depend on $\theta_0$) such that
\begin{equation}
|\partial^3 \log f (x;\theta)/\partial \theta^3|\leq M(x),\ \mbox{for all } x\in A,\ \theta_0-c < \theta < \theta_0 + c
\end{equation}
and
${\mathbb E}_{\theta_0}\left(M(X)\right) < \infty$.
\end{itemize}
Obviously, in  our case $\theta=L_0$, $\Omega=(0,\infty)$ and $f(x;\theta)$ is given by (\ref{param1}). So conditions (i), (ii) and (iii) are fulfilled. On the other hand, (iv) is also fulfilled since the function $f(x;L_0)$ in the integrand has three continuous derivatives with respect to $L_0$.

The validity of condition (v) follows from the direct calculation of the Fisher information quantity which yields 
\begin{eqnarray*}
I(L_0) &=& -{\mathbb E}\left[\ell''(L_0;D) \right]={\mathbb E}\left[-\frac{1}{(L_0+\pi R)^2} + \frac{1}{(L_0 + 2\pi D)^2} \right]\\
&=& \int_0^R \left(-\frac{1}{(L_0+\pi R)^2} + \frac{1}{(L_0 + 2\pi r)^2}\right) \left(\frac{L_0 + 2\pi r}{L_0 R + \pi R^2}\right) dr\\
&=& \frac{1}{L_0 + \pi R}\left[ \frac{1}{2\pi R}\log\left(1+\frac{ 2\pi R}{L_0}\right) - \frac{1}{L_0 + \pi R} \right].
\end{eqnarray*}

As for condition (vi) let us note that
$$
\frac{\partial^3}{\partial L^3} \log (f(x;L))=\frac{2\pi(R-2x)[4\pi^2x^2+2\pi x(3L+\pi R)+3L^2+3\pi LR+\pi^2R^2]}{(L+\pi   R)^3(2\pi x+L)^3}.
$$
 Now, a function $M(x)$ fulfilling condition $(vi)$ in a neighborhood $(L_0-c,L_0+c)$ of $L=L_0$ is, for example,
{\small $$
M(x)=\frac{2\pi (R-2x)\big[4\pi^2x^2+2\pi x\big(3(L_0+c)+\pi R\big)+3(L_0+c)^2+3\pi (L_0+c)R+\pi^2R^2\big]}{(L_0-c+\pi   R)^3(2\pi x+L_0-c)^3},
$$}
which clearly satisfies ${\mathbb E}_{L_0}\left(M(X)\right) < \infty$.

Finally, as a consequence of the asymptotic normality (and asymptotic efficiency) of the maximum likelihood estimator [Theorem 3.10 in Lehmann and Casella, p. 449] we can conclude
\[
\sqrt{n}(\widehat{L}_0-L_0) \stackrel{w}\longrightarrow\mbox{N} (0,I(L_0)^{-1}).
\]

\

\noindent\it The three-dimensional case\rm

We first establish the basic model to be considered in the inference. This is done in the following result, which is the analog of
Proposition \ref{proposition:pv2} for the three-dimensional case. Again, we will provide two alternative expressions for the density $f$ of the random variable $D$, the distance to $S$ from a random uniformly chosen on $B(S,R)\setminus S$. The set $S$ is assumed to belong to the class ${\mathcal P}{\mathcal V}(R,3)$ of compact sets in ${\mathbb R}^3$ with polynomial volume given by
\begin{equation}
V(r;S)=\mu(S)+L_0r+M r^2+4\pi r^3/3,\ 0\leq r<R.\label{pv-3}
\end{equation}

\begin{proposition}
\label{proposition:pv3}
The
 above defined random variable ``distance to the boundary'', $D$, is absolutely continuous with density function
\begin{equation}\label{dens3d}
f(r;L_0,M)=\frac{L_0+2Mr+4\pi r^2}{L_0R+MR^2+\frac{4}{3}\pi R^3},\ 0\leq r<R.
\end{equation}

This density can be alternatively expressed as
\begin{equation}
f(r)=\lambda_1f_1(r)+\lambda_2f_2(r)+\lambda_3f_3(r),\ 0\leq r<R.\label{denspv-3},
\end{equation}
where
{\small $$
\lambda_1=\frac{L_0}{L_0+MR+4\pi R^2/3},\ \lambda_2=\frac{MR}{L_0+MR+4\pi R^2/3},\ \lambda_3=\frac{4\pi R^2/3}{L_0+MR+4\pi R^2/3},
$$}
and, for $i=1,2,3$,   $f_i$ is the density function of a random variable $RZ_i$, where $Z_1$ is uniform on $(0,1)$, $Z_2$ has a distribution $\mbox{Beta}\,(2,1)$ and $Z_3$ is $\mbox{Beta}\,(3,1)$.
\end{proposition}

\noindent \sc Proof\rm: The expression (\ref{dens3d}) follows directly for Proposition \ref{stacho} and the expression (\ref{pv-3}) of the volume function for the sets in ${\mathcal P}{\mathcal V}(R,3)$. The expression (\ref{denspv-3}) is just a reformulation of (\ref{dens3d}).

\

Again, expression (\ref{denspv-3}) can be interpreted in geometric terms: if we think, to fix ideas, that $S$ is a polyhedron, then
$f_1$, $f_2$ and $f_3$ would represent, respectively, the densities of the distances of those points whose projections are inside a face, on an edge and on a vertex.

\

Now, the main results concerning the moment estimators of $L_0$ and $M$ are summarized in the following statement.

\begin{theorem}\label{3d}
Under the assumption (\ref{pv-3}), we have that the estimators of $L_0$ and $M$ by the method of moments from a sample $D_1,\ldots,D_n$ of the distance variable $D$ are
$$
\tilde{L}_0=\frac{2\pi R^ 2}{5}\left[\frac{3R^ 2-12\overline{D}R+10\overline{D^2}}{R^2-6\overline{D}R+6\overline{D^ 2}}\right].
$$
$$
\tilde{M}=-\frac{4\pi R}{5}\left[\frac{3R^2-16\overline{D}R+15\overline{D^2}}{R^2-6\overline{D}R+6\overline{D^2}}\right],
$$
where $\overline{D}$ and $\overline{D^ 2}$ are the sample means of $D_1,\ldots,D_n$ and $D_1^2,\ldots,D_n^2$, respectively. Moreover, if we denote
\begin{equation*}
g_1(u,v)= \frac{2\pi R^2}{5}\left[\frac{3 R^2-12 uR+10v}{R^2-6uR+6v}\right]
\end{equation*}
and
\begin{equation*}
g_2(u,v)=-\frac{4\pi R}{5}\left[\frac{3R^2-16uR+15v}{R^2-6uR+6v}\right]
\end{equation*}
then,
\begin{eqnarray}
\sqrt{n}(\tilde{L}_0-L_0)\stackrel{w}{\longrightarrow} N(0,\sigma^2_{\tilde{L}_0})\mbox{ and }
\sqrt{n}(\tilde{M}-M_0)\stackrel{w}{\longrightarrow} N(0,\sigma^2_{\tilde{M}}),\label{momentos-asint}
\end{eqnarray}
with $\sigma_{\tilde{L}_0}^2= \nabla g_1^t \Sigma_{D,D^2} \nabla g_1$ and
$\sigma_{\tilde{M}}^2= \nabla g_2^t \Sigma_{D,D^2} \nabla g_2$, where
$\Sigma_{D,D^2}$ is the covariance matrix of the vector $(D,D^2)$. The elements of $\Sigma_{D,D^2}$ are
{\small $$Var(D)= \frac{R^2(12\pi^2  R^4+24\pi M R^3+10 M^2 R^2+44\pi L_0 R^2+30 L_0 MR+15 L_0^2)}{20(4\pi R^2+3MR+3L_0)^2},$$}
{\small $$Var(D^2)=\frac{R^4 ( 768\pi^2 R^{4}+1360\pi M R^3+525 M^2 R^2+1920\pi L_0 R^2+1260 L_0 M R+560 L_0^2)}{700( 4\pi R^2+3MR+3L_0)^2},$$}
and
{\small $$Cov(D,D^2)= \frac{ R^3( 16\pi^2 R^4+30\pi M R^3+12 M^2 R^2+48\pi L_0 R^2+32 L_0 M R+15 L_0^2 )}{20( 4\pi R^2+3 M R+3 L_0)^2}.$$}
\end{theorem}

\noindent \sc Proof\rm:
Some elementary calculations lead to
\begin{equation*}
\label{e1}
{\mathbb E}(D)=\int_{0}^ R r\frac{L_0+2Mr+4\pi r^2}{L_0R+MR^2+\frac{4}{3}\pi R^3}dr=\frac{3L_0R+4MR^2+6\pi R^3}{6(L_0+MR+\frac{4}{3}\pi R^2)},
\end{equation*}

\begin{equation*}
\label{e2}
{\mathbb E}(D^2)=\int_{0}^ R r^2\frac{L_0+2Mr+4\pi r^2}{L_0R+MR^2+\frac{4}{3}\pi R^3}dr=\frac{10L_0R^2+15MR^3+24\pi R^4}{30(L_0+MR+\frac{4}{3}\pi R^2)},
\end{equation*}
The estimators $\tilde{L}_0$ and $\tilde{M}$ are then obtained as the solutions of the system of equations
$\overline{D}={\mathbb E}(D)$, $\overline{D^2}={\mathbb E}(D^2)$.

\

With the notation introduced for $g_1$, we have
$$\sqrt{n}(\tilde{L}_0-L_0)=\sqrt{n}\Big(g_1\big(\overline{D},\overline{D^2}\big)-g_1\big({\mathbb E}(D),{\mathbb E}(D^ 2)\big)\Big).$$
Performing a Taylor expansion for $g_1$ at $({\mathbb E}(D),{\mathbb E}(D^2))$ and denoting $v=\big((\overline{D},\overline{D^2})-({\mathbb E}(D),{\mathbb E}(D^2)\big)$, we obtain
$$g_1(\overline{D},\overline{D^ 2})-g_1({\mathbb E}(D),{\mathbb E}(D^2))=\nabla  \big(g_1({\mathbb E}(D),{\mathbb E}(D^2))\big)^t v+r(v).$$

We only need to show $\sqrt{n} r(v)\stackrel{P}{\longrightarrow} 0$. This follows from the fact that $g$ is a function of differentiability class two in a neighborhood of $({\mathbb E}(D),{\mathbb E}(D^2))$. Indeed, we have

$$\frac{\partial g_1}{\partial u}=\frac{12\pi R^3}{5}\frac{R^2-2v}{(R^2-6uR+6v)^2},\
\frac{\partial g_1}{\partial v}=\frac{8 \pi R^3}{5}\frac{3u-2R}{(R^2-6uR+6v)^2}$$

$$\frac{\partial g_1}{\partial u}\big({\mathbb E}(D),{\mathbb E}(D^2)\big)= \frac{3(5L_0-4\pi R^2)(4\pi R^2+3MR+3L_0)}{\pi R^3}.$$

$$\frac{\partial g_1}{\partial v}\big({\mathbb E}(D),{\mathbb E}(D^2)\big)= \frac{5(2\pi R^2-3L_0)( 4\pi R^2+3MR+3L_0)}{\pi R^4}.$$

$$\frac{\partial^2 g_1}{\partial u^2}=\frac{144\pi R^4}{5}\frac{R^2-2v}{\left(R^2-6Ru+6v\right)
^3},
\frac{\partial^2 g_1}{\partial u\partial v}= \frac{24 \pi R^3}{5}\frac{-7R^2+6Ru+6v}{(R^2-6Ru+6v)^3},$$

$$\frac{\partial^2 g_1}{\partial v^2}=\frac{96\pi  R^3}{5}\frac{2R-3u}{(R^2-6Ru+6v)^3}$$

To check the continuity of the second-order derivatives we only have to see that the denominators are not null at  $({\mathbb E}(D),{\mathbb E}(D^2))$. This follows by observing that if we replace $u$ with ${\mathbb E}(D)$ and $v$ with ${\mathbb E}(D^2)$ in $R^2-6Ru+6v$ we get
$$\frac{2\pi R^4}{5}\frac{1}{3L_0+3MR+4\pi R^2},$$
which is not null for $R>0$.

\

To get the asymptotic distribution of $\tilde{M}$ we note that $\tilde{M}= g_2(\overline{D},\overline{D^2})$, where
$$g_2(u,v)=-\frac{4\pi R}{5}\left[\frac{3R^2-16uR+15v}{R^2-6uR+6v}\right].$$
Then
$$\sqrt{n}(\tilde{M}-M)=\sqrt{n}\Big(g_2\big(\overline{D},\overline{D^2}\big)-g_2\big({\mathbb E}(D),{\mathbb E}(D^ 2)\big)\Big).$$
We now make a similar reasoning to that of $\tilde{L}_0$ which requires to calculate the derivatives of first and second order of $g_2$ and to check their continuity at a neighborhood of  $({\mathbb E}(D),{\mathbb E}(D^2))$. This easily follows from
$$\frac{\partial g_2}{\partial u}= -\frac{8\pi R^2(R^2-3v)}{5(R^2-6uR+6v)^2},$$
$$\frac{\partial g_2}{\partial v}= \frac{12\pi R^2(R-2u)}{5(R^2-6uR+6v)^2},$$
$$\frac{\partial g_2}{\partial u}\big({\mathbb E}(D),{\mathbb E}(D^2)\big)= \frac{\left( 32\pi R+15M\right)\left( 4\pi R^2+3MR+3L_0\right) }{\pi  R^{3}},$$
$$\frac{\partial g_2}{\partial v}\big({\mathbb E}(D), {\mathbb E}(D^2)\big) = -\frac{15(2\pi R+M)( 4\pi R^2+3MR+3L_0)}{\pi R^4}.$$
\hfill $\square$

\

We will omit the analysis of the maximum likelihood estimators since the required conditions to ensure asymptotic normality (and asymptotic efficiency) are extremely complicated to check in this case.

\section{An estimation pathology and how to handle it}\label{path}

Our first result in this section applies to the estimator of $L_0$ by the method of moments for the case $d=2$. However, the following discussion suggests that a similar behavior is also present in the other considered cases. The point is that
while it has a quite simple explicit expression and it is asymptotically normal with a explicitly known variance, it has an infinite mean. This could be  seen as  a sort of intrinsic, extreme case of non-robustness. Of course, the problem lies with the samples whose sample mean  is close to the value $R/2$ where the expression of the estimator $\tilde L_0$ goes to infinity. The following result shows that such ``natural outlying'' samples are probable enough to give an infinite expectation for the estimator.

\begin{proposition}\label{L0inf}
In the case $d=2$ the estimator $\tilde{L}_0$ has an infinite expectation.
\end{proposition}
\noindent {\sc Proof:}
Let $X_1,\ldots,X_n$ be i.i.d. random variables taking values on an interval $[a,b]$ and whose density $f$ satisfies $f(u)\geq c >0 $ for all $u\in [a,b]$. For any  function $g(X_1,\ldots,X_n)\geq 0$ we have
\begin{align*}
\mathbb{E}g(X_1,\ldots,X_n) &= \int_{[a,b]^n} g(x_1,\ldots,x_n)\prod_{i=1}^n f(x_i) \, dx_1\cdots dx_n \\
&\geq c^n(b-a)^n  \int_{[a,b]^n} g(x_1,\ldots,x_n)\frac{1}{(b-a)^n}\, dx_1\cdots dx_n\\
&= c^n(b-a)^n\mathbb{E}g(U_1,\ldots, U_n),
\end{align*}
where $U_1,\ldots, U_n$ are  i.i.d. random variables, uniformly distributed on $[a,b]$.
Since the density of the distances satisfies
\[
f(r) \geq \frac{L_0}{L_0R+\pi R^2} >0, \ \mbox{for all}\ r\in [0,R],
\]
we can apply the observation above to deduce the following lower bound:
\[
\mathbb{E} |\tilde{L}_0| = \frac{2\pi R}{3} \mathbb{E} \left|\frac{2R-3\bar{D}}{2\bar{D}-R}\right| \geq
 \frac{2\pi R}{3} \left(\frac{L_0 R}{L_0R+\pi R^2} \right)^n\mathbb{E} \left|\frac{2R-3\bar{U}}{2\bar{U}-R}\right|,
\]
where $\bar{U}=n^{-1}(U_1+\cdots + U_n)$ and $U_1,\ldots, U_n$ are i.i.d. random variables, uniformly distributed on $[0,R]$.  Now, the following equality is easy to check:
\[
\mathbb{E} \left|\frac{2R-3\bar{U}}{2\bar{U}-R}\right| =
\mathbb{E} \left|-\frac{3}{2} +  \frac{R}{4}\, \frac{1}{\bar{U}-R/2}\right|.
\]
Therefore, $\mathbb{E} |\tilde{L}_0| = \infty$ follows as a corollary of the following lemma (applied to $U_1,\ldots, U_n$).
\hfill $\square$

\

\begin{lemma}\label{L0inf1}
Let $X_1,\ldots,X_n$ be i.i.d. random variables taking values on an interval $[a,b]$. Assume their distribution $F$ has a density $f$ such that there exists a  constant $c>0$ with $f(u)\geq c$, for all $u\in [a,b]$. Then,
\[
\mathbb{E} \left| \frac{1}{\bar{X} - \mu}\right| = \infty,
\]
where $\bar{X}=n^{-1}(X_1+\cdots + X_n)$ and $\mu = \mathbb{E}(X_1)$.
\end{lemma}

\noindent {\sc Proof:}
Since,
\[
|\bar{X} - \mu  | = \left| \frac{X_1}{n} - \frac{\mu}{n} + \cdots +
\frac{X_n}{n} - \frac{\mu}{n} \right| \leq n^{-1} \sum_{i=1}^n |X_i - \mu|,
\]
we have
\begin{equation}
\label{eq.lemma}
\mathbb{E} \left| \frac{1}{\bar{X} - \mu}\right| \geq  n\mathbb{E} \left( \frac{1}{\sum_{i=1}^n |X_i - \mu|}\right) =
n\int_0^\infty \mathbb{P}\{\sum_{i=1}^n |X_i - \mu| < 1/t \}\, dt.
\end{equation}

Taking into account that
\[
\bigcap_{i=1}^n \left\{ |X_i - \mu| < \frac{1}{nt}\right\} \subset
\left\{\sum_{i=1}^n |X_i - \mu| < \frac{1}{t}\right\}
\]
we obtain
\begin{eqnarray*}
&&\int_0^\infty \mathbb{P}\{\sum_{i=1}^n |X_i - \mu| < 1/t \}\, dt \geq
 \int_0^\infty \left(\mathbb{P}\{ |X_1 - \mu| < \frac{1}{nt}\} \right)^n \, dt\\ =
&& \int_0^\infty  \left(\int_{\mu-(nt)^{-1}}^{\mu+(nt)^{-1}} f(u)\, du   \right)^n\, dt
\end{eqnarray*}
Using the assumption on the density,
\[
\int_{\mu-(nt)^{-1}}^{\mu+(nt)^{-1}} f(u)\, du  \geq \frac{2c}{nt}.
\]
Therefore,
\[
\int_0^\infty \mathbb{P}\{\sum_{i=1}^n |X_i - \mu| < 1/t \}\, dt \geq
\frac{2^n c^n}{n^n} \int_0^\infty \frac{1}{t^n}\, dt = \infty.
\]

The result follows from this fact, together with (\ref{eq.lemma}). \hfill $\square$

\

A similar conclusion should hold for the moment estimator of the boundary measure in the case $d=3$. As for the maximum likelihood estimators, the analysis is more involved, given the lack of  explicit expressions for such estimators. 

\

\noindent \it Some practical consequences\rm

In summary, we are faced with the following somewhat unusual, interesting situation: as a consequence of the results in Sections 3 and 4, we have some standard, relatively easy to find, estimators $T_n$ which, in spite of being consistent ($T_n\stackrel{P}\rightarrow \theta$) and asymptotically normal ($\sqrt{n}(T_n-\theta)\stackrel{w}\rightarrow N(0,\sigma(\theta))$),  have an infinite expected value. In other words,  as the sample size increases these estimators converge (with an approximately normal distribution) to the true value of the parameter but, still, their estimation error, as measured with the usual $L_1$ or $L_2$ criteria, is infinity. Note that there is no contradiction in that since the weak convergence to the normal distribution, as established by the standard asymptotic normality results, does not entail the corresponding convergence for the moments.  The obvious question is: how to deal with this situation? We have two complementary answers:
\begin{itemize}
\item[(a)] \bf To use an alternative error criterion\rm: It is clear that in this case the usual error criteria for an estimator $T_n$ of a parameter $\theta$ (i.e. ${\mathbb E}|T_n-\theta|$ and ${\mathbb E}(T_n-\theta)^2$), are unsuitable, in the sense that they do not reflect the way in which the estimator $T_n$ approaches the value of the target parameter $\theta$. Then, a possible quite natural alternative would be to use a bounded  error criteria
    \begin{equation}
    d_{BE}(T_n;\theta)={\mathbb E}\left(\frac{|T_n-\theta|}{|T_n-\theta|+1}\right).\label{BE}
    \end{equation}
    The motivation for such an error is very simple as a consequence of the following well-known characterization of the convergence in probability: if $Z_n$, (for $n\in{\mathbb N}$) and $Z$ are random variables, we have
    $$
    Z_n\stackrel{P}\rightarrow Z\ \mbox{if and only if }{\mathbb E}\left(\frac{|Z_n-Z|}{|Z_n-Z|+1}\right)\rightarrow 0,
    $$
This equivalence follows directly from the Dominated Convergence Theorem and the following inequality (combined with Markov's inequality)
$$
{\mathbb P}\{|Z_n-Z|>\epsilon\}\leq {\mathbb P}\{\frac{|Z_n-Z|}{|Z_n-Z|+1}>\frac{\epsilon}{\epsilon+1}\}
$$
\item[(b)] \bf To define suitably modified estimators\rm,  aimed to correct the infinite expectation problem. An idea in this line (for the two-dimensional case $d=2$) is as follows: given an estimator $\hat\lambda$ of the mixture parameter $\lambda$ in (\ref{mixture}),  a natural estimator of $L_0$ is
\begin{equation}
\hat L_0 = \pi R \frac{\hat{\lambda}}{1-\hat{\lambda}} = \pi R \sum_{k=1}^\infty \hat{\lambda}^k.\label{trunc}
\end{equation}
For instance, it is very easy to check that $\mu = \mathbb{E}(D) =  2R/3 - \lambda R/6$. Hence, the estimator of $\lambda$ by the method of moments is $\tilde{\lambda}= 4 - 6\bar{D}/R$.

By Monotone Convergence Theorem, $\mathbb{E}(\hat L_0)<\infty$ if and only if the series $\sum_{k=1}^\infty \mathbb{E}(\hat{\lambda}^k)$ is convergent.
However, even in the case when the series diverges,  it is possible to define an estimator of $L_0$ with finite expectation, although biased, through an appropriate truncation from (\ref{trunc}):
\begin{equation}
\hat{L}_0 = \pi R \sum_{k=1}^K \hat{\lambda}^k.\label{trunc1}
\end{equation}

We may also develop similar ideas \it in the three dimensional case\rm, for the purpose of correcting the finite expectation problem. In the case $d=3$ we have three parameters $\lambda_1$, $\lambda_2$ and $\lambda_3$ defined in (\ref{denspv-3}), which can be estimated by the method of moments as the solutions of the system:
\begin{equation} \label{systlambda}
\begin{array}{cccccc}
\mathbb{E}\big(D\big) = & \frac{R}{2} \lambda_1 &+ &  \frac{2R}{3} \lambda_2&+& \frac{3R}{4}\lambda_3 \\
\mathbb{E}\big(D^2\big) = & \frac{R^2}{3} \lambda_1&+& \frac{R^2}{2}\lambda_2 &+& \frac{3R^2}{5}\lambda_3\\
\mathbb{E}\big(D^3\big)= & \frac{R^3}{4} \lambda_1 &+&\frac{2R^3}{5} \lambda_2&+&\frac{R^3}{2}\lambda_3\\
\end{array}
\end{equation}
If we solve (\ref{systlambda}) we obtain
$$\hat{\lambda}_1=\frac{12(6R^2\overline{D}-20R\overline{D^2}+15\overline{D^3})}{R^3}$$
$$\hat{\lambda}_2=-\frac{30(4R^2\overline{D}-15R\overline{D^2}+12\overline{D^3})}{R^3}$$
$$\hat{\lambda}_3=\frac{20(3R^2\overline{D}-12R\overline{D^2}+10\overline{D^3})}{R^3}$$
Thus, the expressions for the estimators based on the method of moments are
$$\tilde{L}_0=\frac{4\pi R^2}{3} \frac{\hat\lambda_1}{1-\hat\lambda_1}\Big(\frac{1}{1-\hat\lambda_2}\Big)\Big(1-\frac{\hat\lambda_1}{1-\hat\lambda_1}\frac{\hat\lambda_2}{1-\hat\lambda_2}\Big)^{-1}$$
and
$$\tilde{M}=\frac{4\pi R}{3} \frac{\hat\lambda_2}{(1-\hat\lambda_2)}\Big[\frac{1}{1-\hat\lambda_2}\frac{\hat\lambda_1}{1-\hat\lambda_1}\Big(1-\frac{\hat\lambda_1}{1-\hat\lambda_1}\frac{\hat\lambda_2}{(1-\hat\lambda_2)}\Big)^{-1}+1\Big].$$
The truncated versions are
$$\tilde{L}_0=\frac{4\pi R^2}{3 \hat\lambda_2} \sum_{j=1}^K \left[\frac{\hat\lambda_1}{1-\hat\lambda_1}\frac{\hat\lambda_2}{1-\hat\lambda_2}\right]^j$$
and
$$\tilde{M}= \frac{4\pi R}{3}\frac{1}{(1-\hat\lambda_2)}\Big[\sum_{j=1}^K \left[\frac{\hat\lambda_1}{1-\hat\lambda_1}\frac{\hat\lambda_2}{1-\hat\lambda_2}\right]^j+\hat\lambda_2\Big].$$

The practical use of these estimators will require some study on the optimal values of $R$ and $K$. This question will not be considered here.

\end{itemize}

\

In the simulation results of the next section we will incorporate the ideas (a) and (b): the performance of the different estimators (moments and maximum likelihood) and   that of their ``truncated'' versions (\ref{trunc1})  have been checked using the error criterion (\ref{BE}).

\

\noindent \it Some further consequences of the mixture representation (\ref{mixture})\rm

As an additional advantage of (\ref{mixture}),   the maximum likelihood estimator of $\lambda$ can be easily computed using the EM-algorithm:

\begin{enumerate}
\item Initial step: $\hat{\lambda}^{(0)} = 0.5$.

\item Iterate until convergence:

\begin{enumerate}

\item E-step. For $i=1,\ldots,n$, let $Y_i$ be the (unobservable) random variable which indicate if $D_i$ has been drawn from $f_1$ or $f_2$ ($Y_1=1$ and $Y_i=0$, respectively). Compute, using Bayes formula,
\[
Y_{i,k} = \mathbb{E}(Y_i | D_i, \hat{\lambda}^{(k)}) = \frac{\hat{\lambda}^{(k)} f_1(D_i)}{\hat{\lambda}^{(k)} f_1(D_i) +
(1-\hat{\lambda}^{(k)}) f_2(D_i)},
\]
and define
\[
Q(\lambda, \hat{\lambda}^{(k)}) = \sum_{i=1}^n [Y_{i,k} \log \lambda
+ (1- Y_{i,k})\log (1-\lambda)].
\]

\item M-step. Find the value $\hat{\lambda}^{(k+1)}$ that maximizes $Q(\cdot,\hat{\lambda}^{(k)})$. It is straightforward to show that $\hat{\lambda}^{(k+1)} = n^{-1}\sum_{i=1}^n Y_{i,k}$.

\end{enumerate}
\end{enumerate}

Notice that if we could observe the variables $Y_i$, the maximum likelihood estimator of $\lambda$ would be $n^{-1}\sum_{i=1}^n Y_i$. Each step of the algorithm uses essentially this formula but replacing $Y_i$ with the corresponding expected value given the current value of $\lambda$.

\section{Simulation results}\label{sim}
\normalsize

\

\noindent \it The two-dimensional case\rm

We have carried out a small simulation to illustrate some aspects of the behavior of the estimators defined in the previous section. Consider the set $S$ defined as the union of two disjoint circles with centers at $(-2.75,0)$ and $(2.75,0)$, and common radius equal to $0.25$. These values imply that the reach of $S$ is 2.5 and its perimeter is $L_0=\pi$.

We have compared four estimators of $L_0$: moments estimator, maximum likelihood and the ``truncated'' versions of them defined in (\ref{trunc1}). The outputs in Tables 1-3 below are based on $B=2000$ replications.

\begin{table}[ht]
\begin{center}
\begin{tabular}{|r|rrrr|rrrr|}
  \hline
  & $R=1$ & $R=1$ & $R=1$ & $R=1$ & $R=2$ & $R=2$ & $R=2$ & $R=2$ \\
  \hline
$n$ & MLE & TMLE & MOM & TMOM & MLE & TMLE & MOM & TMOM \\
  \hline
  100 & 0.530 & 0.505 & 0.538 & 0.514 & 0.536 & 0.531 & 0.554 & 0.549 \\
  300 & 0.411 & 0.389 & 0.420 & 0.399 & 0.421 & 0.417 & 0.443 & 0.438 \\
  500 & 0.364 & 0.343 & 0.371 & 0.351 & 0.374 & 0.370 & 0.390 & 0.386 \\
  700 & 0.335 & 0.316 & 0.345 & 0.325 & 0.334 & 0.330 & 0.352 & 0.349 \\
  1000 & 0.308 & 0.289 & 0.314 & 0.294 & 0.298 & 0.295 & 0.317 & 0.314 \\
  20000 & 0.097 & 0.095 & 0.100 & 0.099 & 0.098 & 0.097 & 0.106 & 0.106 \\
   \hline
\end{tabular}
\caption{Error ($d_{BE}$) averages  over 2000 replications for the maximum likelihood estimator (MLE), the moments estimator (MOM) and their respective truncated versions (TMLE and TMOM). Truncated versions correspond to $K=5$ ($n\neq 20000$) and $K=8$ ($n=20000$). The value of the parameter is $L_0=\pi$. In the E-M algorithm we have taken $10^{-5}$ as the tolerance for $\lambda$.}
\label{d21}
\end{center}
\end{table}

\begin{table}[ht]
\begin{center}
\begin{tabular}{|r|rrrr|rrrr|}
  \hline
  & $R=1$ & $R=1$ & $R=1$ & $R=1$ & $R=2$ & $R=2$ & $R=2$ & $R=2$ \\
  \hline
$n$ & MLE & TMLE & MOM & TMOM & MLE & TMLE & MOM & TMOM \\
  \hline
  100 & 3.07 & 2.98 & 3.07 & 2.99 & 3.16 & 3.14 & 3.21 & 3.20 \\
  300 & 3.10 & 3.01 & 3.09 & 3.00 & 3.09 & 3.08 & 3.09 & 3.07 \\
  500 & 3.14 & 3.04 & 3.14 & 3.05 & 3.11 & 3.10 & 3.15 & 3.14 \\
  700 & 3.10 & 3.00 & 3.12 & 3.02 & 3.15 & 3.14 & 3.13 & 3.12 \\
  1000 & 3.14 & 3.04 & 3.14 & 3.04 & 3.11 & 3.10 & 3.11 & 3.10 \\
  20000 & 3.14 & 3.13 & 3.15 & 3.13 & 3.14 & 3.14 & 3.14 & 3.14 \\
   \hline
\end{tabular}
\caption{Medians over 2000 replications for the maximum likelihood estimator (MLE), the moments estimator (MOM) and their respective truncated versions (TMLE and TMOM). The value of the parameter is $L_0=\pi$.}
\label{d22}
\end{center}
\end{table}

\begin{table}[ht]
\begin{center}
\begin{tabular}{|r|rrrr|rrrr|}
  \hline
  & $R=1$ & $R=1$ & $R=1$ & $R=1$ & $R=2$ & $R=2$ & $R=2$ & $R=2$ \\
  \hline
$n$ & MLE & TMLE & MOM & TMOM & MLE & TMLE & MOM & TMOM \\
  \hline
  100 & 1.85 & 1.67 & 1.92 & 1.77 & 2.03 & 2.00 & 2.16 & 2.11 \\
  300 & 1.12 & 1.00 & 1.15 & 1.03 & 1.16 & 1.14 & 1.30 & 1.28 \\
  500 & 0.89 & 0.79 & 0.93 & 0.83 & 0.94 & 0.92 & 1.02 & 1.00 \\
  700 & 0.75 & 0.68 & 0.80 & 0.71 & 0.76 & 0.75 & 0.84 & 0.82 \\
  1000 & 0.66 & 0.59 & 0.68 & 0.61 & 0.64 & 0.63 & 0.71 & 0.70 \\
  20000 & 0.14 & 0.14 & 0.14 & 0.14 & 0.14 & 0.14 & 0.16 & 0.16 \\
   \hline
\end{tabular}
\caption{Median absolute deviations over 2000 replications for the maximum likelihood estimator (MLE), the moments estimator (MOM) and their respective truncated versions (TMLE and TMOM). The value of the parameter is $L_0=\pi$.}
\label{d23}
\end{center}
\end{table}

In order to properly interpret these outputs, we should keep in mind that the atypical behavior of our estimators requires to modify the usual approach in most simulation studies. In particular, the average of the estimated values along the 2000 runs is no longer here a representative value of the estimator's performance when the corresponding theoretical value of the expectation is infinity. In those cases, the empirical average wouldn't show any apparent improvement as the sample size increases, in spite of the fact that the estimator does converge to the true value of $L_0$. Then, in Table \ref{d21} we just consider the bounded error measure $d_{BE}$ defined in (\ref{BE}): the values of this error measure improve  as $n$ increases, thus showing
in numerical terms the consistency of the estimators.

Table \ref{d22} gives an idea of the evolution of each estimator but replacing the average value over the 2000 replications with the corresponding median, thus avoiding the infinite expectation problem; recall that the asymptotic normality entails the convergence of the respective medians to the limit median (but not the moment convergence). Finally, Table \ref{d23} gives the median absolute deviation (MAD) for the estimators under study. It is re-scaled in the usual way to get consistency in the Gaussian case.

The choice of the values for $n$ is aimed to show this progressive improvement starting from a small/moderate sample size $n=100$ until the large value $n=20000$. The motivation for this latter choice is to check the ``asymptotic'' performance of our estimation method as a numerical (stochastic) algorithm to approximate $L_0$ even in those cases where $S$ is completely known\/ \rm (though possibly with a complicated shape). In those situations the required samples would be obtained by a Monte Carlo procedure, so that the sample size is just limited by our computational power.

As a further consequence of the atypical situation we have found, let us note that when an estimator fulfills ${\mathbb E}(T_n)=\infty$ and still $\sqrt{n}(T_n-\theta)\stackrel{w}\rightarrow N(0,\sigma(\theta))$, the asymptotic variance must be carefully interpreted just as the variance of the asymptotic distribution. This is not the same as the approximate variance (for $n$ large) of $\sqrt{n}(T_n-\theta)$ (which is again infinity).

\begin{figure}[h]
\begin{center}
\includegraphics[height=7cm]{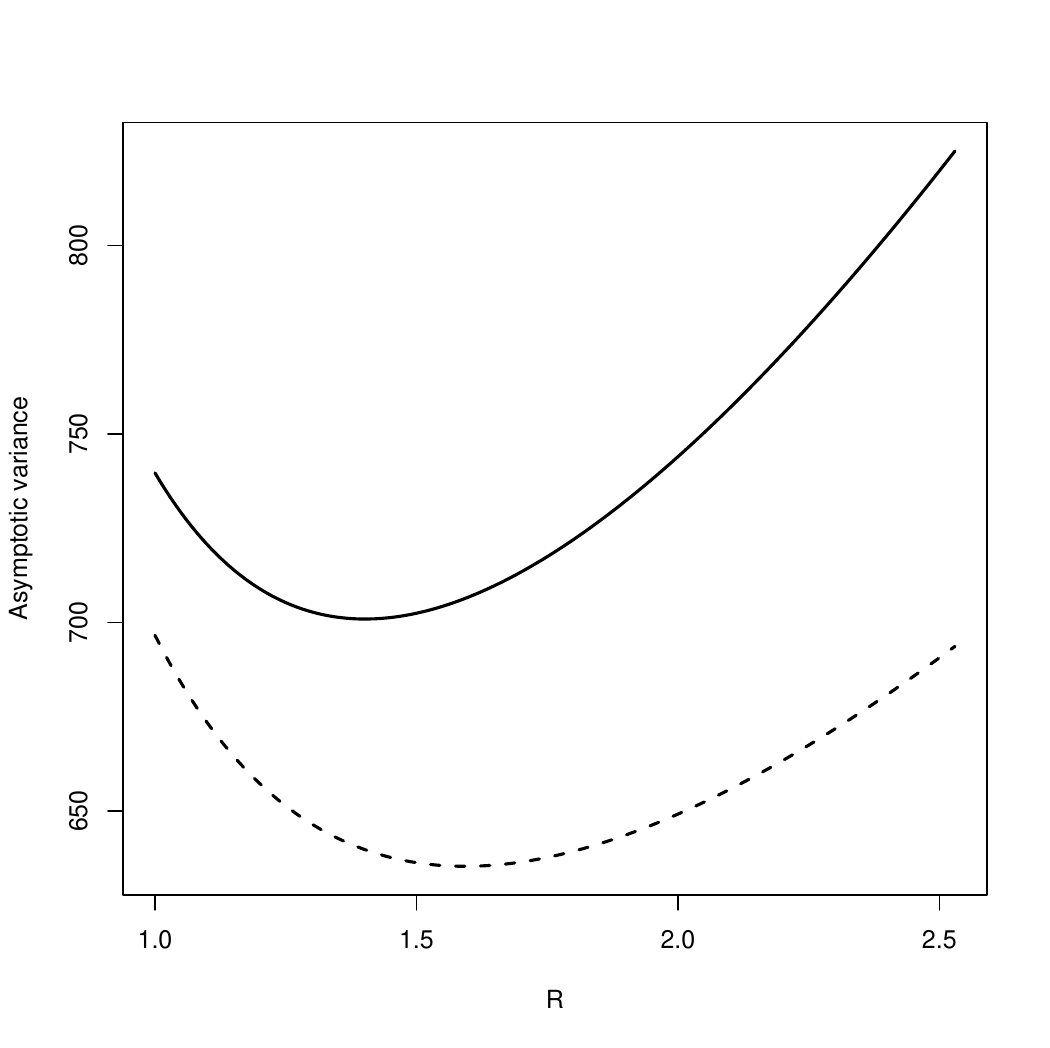}
\end{center}
\caption{Asymptotic standard deviations of the method of moments estimator (solid line) and the maximum likelihood estimator (dashed line) as a function of $R$.}
\label{fig:asympsd}
\end{figure}

In any case, the graphical representation of the asymptotic variances for the estimators obtained by maximum likelihood and the method of moments provides some interesting insights. We have computed
the asymptotic standard deviations  $\sigma_{\tilde{L}_0}$ and $\sigma_{\widehat{L}_0}$, for $R$ ranging between 1 and the reach of $S$ (the set defined at the beginning of this section). The results are displayed in Figure \ref{fig:asympsd}.

It can be seen that
\begin{itemize}
\item[(i)] the asymptotic variance is consistently smaller for the maximum likelihood estimator.
\item[(ii)] The values of the asymptotic variance (in both cases) depend on $R$ in a very natural way, which could be foreseen from the
mixture representation (\ref{mixture}): values of $R$ too small (resp. too large) lead to infra-estimate (resp. over-estimate) the curvature in the boundary of $S$. To see this in the simpler case $S=[0,1]^2$, a large value of $R$ would produce too many points projecting on the vertices of $S$ and a small $R$ would lead to very few points of this type. Then one could say that for each set $S$ (or rather for each volume function) one has an optimal value of $R$.
\item[(iii)] Of course, the truncated estimators will fail to be consistent, unless we would take $K=K_n\to \infty$ in a suitable way. Also, these estimators show a sort of ``bias in the median''  in the sense that their medians over the 2000 replications have often (see Table 2) a larger deviation of the target than the medians of the original (MOM and MLE) estimators.
\end{itemize}


\noindent \it The three-dimensional case\rm

\begin{figure}[h]
\centerline{\includegraphics[scale=.4]{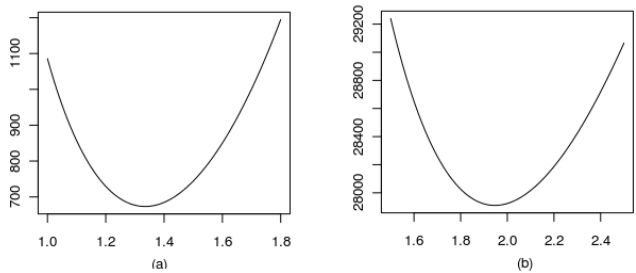}}
\caption{\footnotesize{(a) Asymptotic variance of $\tilde{L}_0$ for a cone of height 1. (b) Asymptotic variance of $\tilde{M}$ for a cone of height 1.}}
\label{moments}
\end{figure}

As in the two-dimensional case, we have computed the asymptotic standard deviations
$\sigma_{\tilde{L}_0}$ and $\sigma_{\tilde{M}}$. The results are displayed in Figure \ref{moments}. They are qualitatively very similar to those for $d=2$: again, the dependence of the results from the value of $R$ illustrates to what extent it is adequate the proportion of points projecting on flat, edgy or corner areas.

\begin{table}[h]
\begin{center}
\footnotesize{
\begin{tabular}{|c|c|c|c|c|}
\hline
                   &  $R=1.3$             & $R=1.9$          &  $R=1.3$         & $R=1.9$       \\
      n            &  MOM $\tilde{L}_0$   & MOM  $\tilde{M}$ &  MLE $\hat{L}_0$ & MLE  $\hat{M}$\\
\hline
5000             & 3.10             & 7.16         & 3.15         & 7.14      \\
20000            & 3.12             & 6.93         & 3.14         & 7.04      \\
40000            & 3.14             & 6.98         & 3.14         & 7.04      \\
\hline
\end{tabular}
\caption{Medians over 2000 replications for the MOM and MLE estimators for $d=3$.}
\label{Mom_med_cono} }
\end{center}
\end{table}

\begin{table}[h]
\begin{center}
\footnotesize{
\begin{tabular}{|c|c|c|c|c|}
\hline
                 &    $R=1.3$          & $R=1.9$             &  $R=1.3$         & $R=1.9$     \\
  n              &  MOM $\tilde{L}_0$  & MOM  $\tilde{M}$    &  MLE $\hat{L}_0$ & MLE  $\hat{M}$ \\
\hline
5000             & 0.65                & 4.25                & 0.38             & 1.89     \\
20000            & 0.32                & 2.06                & 0.19             & 0.96     \\
40000            & 0.23                & 1.52                & 0.13             & 0.71     \\
\hline
\end{tabular}
\caption{ Median absolute deviation over 2000 replications for the MOM and MLE estimators for $d=3$.}
\label{Mom_mad_cono} }
\end{center}
\end{table}

\begin{table}[h]
\begin{center}
\footnotesize{
\begin{tabular}{|c|c|c|c|c|}
\hline
                   &  $R=1.3$                 & $R=1.9$               &  $R=1.3$           & $R=1.9$     \\
n                  &  MOM $\tilde{L}_0$       & MOM  $\tilde{M}$      &  MLE $\hat{L}_0$ & MLE  $\hat{M}$\\
\hline
5000              & 0.31                 & 0.67                  & 0.22           & 0.52    \\
20000              & 0.19                 & 0.54                  & 0.13           & 0.38    \\
40000              & 0.14                 & 0.46                  & 0.10           & 0.32    \\
\hline
\end{tabular}
\caption{Error ($d_{BE}$) average over 2000 replications for the MOM and MLE estimators for $d=3$.}
\label{Mom_error_cono} }
\end{center}
\end{table}

Regarding the numerical experiments, the set $S$ we have considered is a cone with height 1 and angle $\pi/3$ for which $L_0=\pi$ and $M=6.9404$. The results are displayed in Tables \ref{Mom_med_cono}, \ref{Mom_mad_cono} and \ref{Mom_error_cono}. In this case, we just provide  the outputs corresponding to the maximum likelihood (MLE) and the moment-based (MOM) estimators.

In the three-dimensional case the required sample sizes are much larger than those needed for $d=2$. Again, this is just a consequence of the intrinsic difficulty of the considered problem. On the other hand,  the parameters  ($R$, $K$, the tolerance in the
$E$-$M$-algorithm) turn out  to be critical for the estimators performance. A detailed study of the optimal choice of these parameters would require extensive numerical simulation, far beyond the scope of this work.

\newpage

\section{Some final remarks}

The statistical model under study, based on distances to the boundary of the body from randomly chosen outside points is of potential interest in remote sensing in those situations where the distance to the object of interest is measured, for instance, from radar or sonar techniques. Thus, while our results apply in principle to some simplified situations, they could shed some light on some theoretical and  geometric aspects of these methodologies. Let us note that the estimation procedures in the above section only depend on the specific shape of $S$, through the expression of the volume function.

From the methodological point of view, the theoretical development in Sections \ref{stat} and \ref{path} provide a curious example where the classical point estimation theory applies nicely.  Thus, for example, the estimator for the method of moments, together with its asymptotic variance, can be explicitly obtained in closed form. The likelihood function is also explicitly found (and it is not difficult to handle) and Fisher information measure (which provides the asymptotic variance for $d=2$) can be also easily calculated.

However, such an apparently simple approach leads to a somewhat surprising scenario where some natural estimators have an infinite expectation. This entails some non-trivial challenges, as we have discussed in Section \ref{path}.

     The numerical outputs of Section \ref{sim} show that the estimation problem is intrinsically difficult, so that relatively large samples are required. The use of large sample sizes would not represent any major problem in many practical situations where our distance data could be obtained in a simple an inexpensive way. These would include the cases where $S$ is known and the samples are obtained by Monte Carlo simulation. In those situations our procedure might be interpreted as a sort of stochastically-\-based numerical method to approximate the unknown quantities.

It should be also stressed that our results apply to the case that the outside points, from which the distances $D_i$ are calculated,  come from a uniform sample of a crown outside $S$. It is natural to ask to which extent the results rely on the uniformity assumption. Thus, the standard robustness techniques (as found, for instance in the classical book by Huber, 1980) are also in order here.

\section*{Acknowledgements}

We are very grateful to the comments and criticisms from an anonymous referee which led to a substantially improved version of this manuscript.

\section*{References}

\begin{list}{}{\leftmargin .5cm\listparindent -.5cm}
\item\hspace{-.2cm}

\vspace{-.3cm}
{\sc Ambrosio, L., Colesanti, A}. and {\sc Villa, E.} (2008). Outer Minkowski content for some classes of closed sets. \emph{Math. Ann.} \bf 342\rm, 727--748.

{\sc Armend\'ariz, I., Cuevas, A}. and {\sc Fraiman, R.} (2009). Nonparametric estimation of boundary measures and
related functionals: asymptotic results. \emph{Adv. in Appl. Probab.}, \bf 41\rm, 311-322.

{\sc Cuevas, A}. and {\sc Fraiman, R}. (2009). Set estimation. In \it New perspectives in stochastic geometry\rm, W.G. Kendall and I. Molchanov, eds., pp. 366-389. Oxford University Press.

{\sc Cuevas, A., Fraiman, R}. and {\sc Rodr\'{\i}guez-Casal, A}. (2007).
A nonparametric approach to the estimation of lengths and surface areas.
\emph{Ann. Statist.\/},~{\bf 35}\rm, 1031--1051.

{\sc Cuevas, A.}  and {\sc Rodr\'{\i}guez-Casal, A}. (2004). On boundary estimation.  \emph{Adv. in Appl. Probab.}, \bf 36\rm, 340-354

{\sc Delfour, M.C.} and {\sc Zolsio, J.P}. (2001).
\it Shapes and Geometries\rm.
Society for Industrial and Applied Mathematics (SIAM), Philadelphia.

{\sc Dey, T.K}. (2007). \it Curve and Surface Reconstruction\rm. Cambridge University Press.

{\sc Federer, H}. (1959). Curvature measures.
\emph{Trans. Amer. Math. Soc.\/},~{\bf 93}\rm, 418--491.

{\sc Genovese, C.R., Perone-Pacifico, M., Verdinelli, I.} and {\sc Wasserman, L}. (2012a). Minimax manifold estimation. \emph{J. Mach. Learn. Res}. \bf 13\rm,
1263--1291.

{\sc Genovese, C.R., Perone-Pacifico, M., Verdinelli, I.} and {\sc Wasserman, L}. (2012b). Manifold estimation and singular
deconvolution under Hausdorff loss. \emph{Ann. Statist.\/},~{\bf 40}\rm, 941--963.


{\sc Hatcher, A. } (2002). \it Algebraic Topology\rm. Cambridge University Press.

{\sc Heveling., M., Hug, D}. and  {\sc Last, G}. (2004). Does polynomial parallel volume imply convexity? \emph{Math. Ann}., \bf 328\rm, 469-479.

{\sc Hug, D., Last, G}. and {\sc Weil, W}. (2004). A local Steiner-type formula for general closed sets
and applications. \emph{Math. Z}., \bf 246\rm, 237--272.

{\sc Huber, P.J}. (1980). \it Robust Statistics\rm. Wiley, New York.

{\sc Jim\'enez, R}. and {\sc Yukich, J.E}. (2011).
Nonparametric estimation of surface integrals.
\emph{Ann. Statist.\/},~{\bf 39}\rm, 232--260.

{\sc Lehmann, E.L}. and {\sc Casella, G.}. (1998). \it Theory of Point Estimation (2nd edition)\rm. Springer, New York.

{\sc Niyogi, P}., {\sc Smale, S}. and {\sc Weinberger, S}. (2008). A topological view of unsupervised
learning from noisy data. Manuscript available at 
\begin{verbatim}
http://people.cs.uchicago.edu/~niyogi/papersps/noise.pdf
\end{verbatim} \rm 
(last accessed:
12 december, 2012).

{\sc Pateiro-L\'opez, B}. and {\sc Rodr\'{\i}guez-Casal, A}. (2008).
Length and surface area estimation under convexity type
restrictions.
\emph{Adv. in Appl. Probab.\/},~{\bf 40}\rm, 348--358.

{\sc Stach\'o, L.L}. (1976). On the volume function of parallel sets. \it Acta Sci. Math., \bf 38\rm, 365--374.

{\sc Steiner, J.} (1840). \"Uber parallele Fl\"achen. \it Monatsbericht der Akademie der Wissenschaften
zu Berlin\/ \rm pp. 114-118.

{\sc Villa, E}. (2009).
On the outer Minkowski content of sets.
\emph{Ann. Mat. Pura Appl.},~{\bf 188}\rm, 619--630.

{\sc Walther, G.} (1997).
Granulometric smoothing.
\emph{Ann. Statist.} {\bf 25}\rm, 2273--2299.

\end{list}

\end{document}